Manuel DE LA SEN[*]

# STABILITY OF DELAYED SYSTEMS MODELED BY FRACTIONAL MODELS

This paper discusses linear fractional representations (LFR) of parameter-dependent nonlinear systems with real-rational nonlinearities and point-delayed dynamics. Sufficient conditions for robust global asymptotic stability both independent of and dependent on the delays are investigated via linear matrix inequalities. Such inequalities are obtained from the values of the time-derivatives of appropriate Lyapunov functions at all the vertices of the polytope which contains the parametrized uncertainties

## 1  INTRODUCTION

Time-delay systems are very common in nature like, for instance, related to transportation problems, population growing and signal transmission methods  (see, for instance, [1-2] and references therein). The stability and stabilization of those systems have been studied in the literature in connection, for instance, with  Lyapunov  theory. Some of the related results are referred to either as being independent of the sizes of the delays or as dependent of those sizes. Within this last class of results, they merit special attention those related to the characterization of the first interval of admissible delay sizes allowing stabilization. On the other hand, the most involved group of results to obtain is that related to internal delays (i.e. in the state) since its associate dynamics possess infinitely many modes in general.  In this paper, we consider a parameter-dependent (in general, nonlinear and time-varying) system subject to a finite set of point delays which may be, in general, defined by real-rational  nonlinearities, whose parameter vector $H_\infty$ is restricted to lie in a polytope $\Theta \in \mathbf{R}^n$ containing the origin, This is called a so-called polytopic delayed system following the nomenclature used for delay-free systems in [3]]. The results developed in the following might be still applied if the set $\Theta$ is not a polytope after replacing  it by some polytope $\Theta_{poly} \supset \Theta$ . The main arguments used to  develop the formalism are based on the fact that the polytope where the parameters belong to defines affine function matrices of vertices which may be calculated from those of the original polytope $\Theta \in \mathbf{R}^n$ of parametrized uncertainties. In the following, the robust global asymptotic stability of such a polytopic delayed system subject to point delays is investigated via Lyapunov theory.

[*] M. De la Sen, Instituto de Investigación y Desarrollo  de Procesos, Department of Systems Engineering and  Automatic Control, Faculty of Science.  Leioa  (Bizkaia). Aptdo. 644- Bilbao, 48080- Bilbao, Spain, e-mail: wepdepam@lg.ehu.es

I.1. Notation.- $\mathbf{R}^{m \times n}(\mathbf{C}^{m \times n})$ is the set of real (complex) m × n matrices and $P = P^T > 0$ stands for a real symmetric positive-definite matrix.
- For a given set S, one defines $\sigma S = \{\sigma s : s \in S\}$ if $\sigma$ is a positive number.
- The convex hull of complex m × n matrices $(\vartheta_1, \vartheta_2, ... \vartheta_l)$, $\vartheta_i \in \mathbf{C}^{m \times n}$ is

$$\mathbf{Co}\{\vartheta_1, \vartheta_2, ... \vartheta_l\} = \left\{ \vartheta : \vartheta = \sum_{i=1}^{l} \lambda_i \vartheta_i, \sum_{i=1}^{l} \lambda_i \leq 1, \lambda_i \geq 0 \right\}$$

- $I_m$ is the m-identity matrix with the subscript being omitted if its size follows directly from context.
- If $\Theta$ is a polytope containing the origin and $\Delta_i(\theta)$ $(i = \overline{0,r})$, $r \geq 0$ being an integer, are real-valued rational matrix functions of any order of $\theta$ then $\Delta_i = \{\Delta_i(\theta) : \theta \in \Theta\}$ and $\Delta = \Delta_0 \times \Delta_1 \times ... \Delta_r$ are polytopes of $v_i$ vertices $\Delta_i^{(k_i)}$, $k_i = \overline{0, v_i}$ ; $i = \overline{0, r}$ ; and $\Delta^{(k_0, k_1, ..., k_r)} = \Delta_0^{(k_0)} \times ... \times \Delta_r^{(k_r)}$, respectively, where $´\times´$ denotes the Cartesian product of matrices (considered as sets). In our context, $\Theta$ is the polytope where the system parameters belong to while $\Delta_i$ is the polytope where the rational matrix function $A_i(\theta)$, defining the dynamics of the i-th delay $h_i$ ($A_0(\theta)$ describes the delay-free dynamics; i.e. $h_0 = 0$) as the parameter vector $\theta$ runs over $\Theta$ ; $i = \overline{0, r}$.

## 2. LINEAR FRACTIONAL DESCRIPTIONS

Consider the parameter dependent system subject to r point delays $h_i$ $(i=\overline{1,r})$:

$$\dot{x}(t) = \sum_{i=0}^{r} A_i(\theta(t)) x(t - h_i) + B(\theta(t)) u(t) \tag{1.a}$$

$$y(t) = C(\theta(t)) x(t) + D(\theta(t)) u(t) \tag{1.b}$$

where $h_0 = 0$, $x(t) \in \mathbf{R}^n$, $u(t) \in \mathbf{R}^{n_u}$, $y(t) \in \mathbf{R}^{n_y}$ are the state, input and measurable output signals, respectively, and $A_i (i = \overline{0, r})$, B, C and D are real-valued rational functions of time-varying parameter vector $\theta(t)$ with $\theta \in \Theta$ for all $t \geq 0$ with $\theta(t) = (\theta_1(t), \theta_2(t), ..., \theta_m(t))^T$ such that the real vector associated function $\theta : \Theta (\subset \mathbf{R}^m) \times [0, t) \to \mathbf{R}^m$ is defined in such a way that (1) has a mild solution on $[0, t)$ for all $t \geq 0$ for any absolutely continuous function $\varphi : [-h, 0] \to \mathbf{R}^n$ of initial conditions $x(t) \equiv \varphi(t)$, $t \in [-h, 0]$ with $h = \underset{1 \leq i \leq r}{\text{Max}}(h_i)$. One defines:
- The unforced system (1) is robustly globally asymptotically stable if $\|x(t)\|$ is uniformly bounded and $\lim_{t \to \infty} x(t) = 0$ if $u \equiv 0$ for any bounded $x(0)$. The system (1) is robustly stabilizable if there exists an output-feedback realizable control law $u(t) = K(y, \theta(t), t)$ such that the closed-loop system is robustly globally asymptotically stable. For terminology simplicity, since no confusion is expected, we refer in the sequel to robust global asymptotic stability simply as "robust stability".

- The robust stability (stabilizability) margin of (1) for an uncertainty set is $\sigma_m(\rho_m) = \text{Sup}$ $\{\sigma(\rho) : \text{System (1) is robustly stable (stabilizable) over } \gamma\Theta \text{ for all } \gamma \in [0, \sigma]\}$ ($\gamma \in [0, \rho]$). Now, first consider the unforced version of (1) given by,

$$\dot{x}(t) = \sum_{i=0}^{r} A_i(\theta(t))x(t - h_i) \; ; \; y(t) = C(\theta(t))x(t) \tag{2}$$

<u>First LFR</u> : Since $A_i(\theta(t))$ $(i = \overline{0, r})$ is a real-valued rational matrix function of $\theta(t)$, the LFR description of each matrix function $A_i(\theta(t))$ exists for some appropriate matrices $A_{oi}$, $B_{qi}$, $D_{pqi}$ $(i = \overline{0, r})$:

$$A_i(\theta(t)) = A_{0i} + B_{qi}\Delta_i(\theta(t))(I_{d_i} - D_{pqi}\Delta_i(\theta(t)))^{-1}C_{pi} \tag{3}$$

for any $\Delta_i((\theta(t))$ such that the well-posedness condition $\text{Det}(I_{d_i} - D_{pqi}\Delta_i(\theta(t))) \neq 0$, $\forall \theta \in \Theta$, all $t \geq 0$ where $I_{d_i}$ is the $d_i$ identity matrix $(i = \overline{0, r})$. In the following, the explicit dependence of $\theta(t)$ on time is omitted in the notation for the shake of simplicity when no confusion is expected. If (2) is quadratically stable then $A_{0i}$ $(i = \overline{0, r})$ are strictly Hurwitzian (i.e. with all their eigenvalues in Re s < 0). A state-space realization of the state evolution of the dynamical system (2) using (3) is

$$\dot{x}(t) = \sum_{i=0}^{r} \left(A_{0i}x(t - h_i) + B_{qi}q_i(t - h_i)\right)$$
$$p_i(t) = C_{pi}x(t) + D_{pqi}q_i(t) = (I - D_{pqi}\Delta_i(\theta))^{-1}C_{pi}x(t)$$
$$q_i(t) = \Delta_i(\theta)p_i(t) = \Delta_i(\theta)(I - D_{pqi}\Delta_i(\theta))^{-1}C_{pi}x(t)$$
$$\Delta_i(\theta) = \text{Diag}(\theta_1 I_{s_{1i}}, ..., \theta_m I_{s_{mi}}) \tag{4}$$

where $q_i \in \mathbf{R}^{d_i}$; $p_i \in \mathbf{R}^{d_i}$ $(i = \overline{0, r})$ where $s_i = \max_{1 \leq i \leq r}(s_{ki})$ $(i = \overline{0, r})$ and the polytope $\Delta_i = \{\Delta_i(\theta) : \theta \in \Theta\}$ $(i = \overline{0, r})$ are, respectively, the i-th LFR degree with respect to the delay $h_i$ and the polytope of $v_i$ vertices $\Delta_i^{(k)}$ $(k = \overline{1, v_i}; i = \overline{0, r})$ for parametrizations of $A_i(\theta)$ $(i = \overline{0, r})$. In particular, $s_0$ is the LFR degree of the delay-free system (4) (i.e., for the case when $A_i = 0$; $i = \overline{1, r}$) of parametrization within a polytope $\Delta_0$ of $v_0$ vertices $\Delta_0^{(k_0)}$ $(k = \overline{1, v_0})$. All those vertices will become crucial in the subsequent stability analysis in the case of convex problems since it would be sufficient to check stability conditions from matrix inequalities at all the set of distinct vertices. Note that the role of the signals $q_{(.)}$ is that of normalized " equivalent inputs" from the " equivalent outputs" $p_{(.)}$ trough normalization matrices $\Delta_{(.)}$ which depend on the values of the uncertainty parameter vector. The LFR representation (4) will be the main tool of the robust stability analysis and stabilization procedure proposed in this paper. Note also that the LFR (4) of (2) is valid, in particular, for the case of commensurate delays $h_i = i\,h$. $(i = \overline{0, r})$. The uncertainty $\Delta(t - \underline{h}) = \text{Diag}(\Delta_0(t - h_0), ..., \Delta_r(t - h_r))$ is absorbed in the forward loop while the identity operator plays the role of the uncertainty. An alternative single-input LFR to (4) for the state of the system (2) is proposed in the following:

<u>Second LFR</u>: Decompose the rational real-valued matrix function $A(\theta(t)) = \sum_{i=0}^{r} A_i((\theta(t))) U(t - h_i)$, U (t) being the unity Heaviside function, as

$$A(\theta(t)) = \sum_{i=0}^{r} \left( A_{0i} + B_q \Delta(\theta(t))(I - D_{pq}\Delta(\theta(t)))^{-1} C_{pi} U(t - h_i) \right) \quad (5)$$

with $A_{(.)0}$, $B_q$, $D_{pq}$ and $C_{p(.)}$ are real matrices of appropriate sizes. Under well-posedness; i.e. $\operatorname{Det}(I - D_{pq}\Delta(\theta)) \neq 0$, this leads to the single-input multi-output LFR of the state equation of (2):

$$\dot{x}(t) = \sum_{i=0}^{r} A_{0i} x(t-h_i) + B_q q(t) \quad ; \quad q(t) = \Delta(\theta) p(t)$$

$$p(t) = \sum_{i=0}^{r} C_{pi} x(t-h_i) + D_{pq} q(t)$$

$$\Delta(\theta) = \operatorname{Diag}(\theta_1 I_{s_1}, \theta_2 I_{s_2}, ..., \theta_m I_{s_m}) \quad (6)$$

where $q \in \mathbf{R}^d$, $p \in \mathbf{R}^d$, $B_q \in \mathbf{R}^{n \times d}$, $C_{pi} \in \mathbf{R}^{d \times n}$ $(i = \overline{0, r})$, and $\Delta(\theta) \in \mathbf{R}^{d \times d}$ is a matrix function of the time-varying parameter vector $\theta$ (t) in $\Theta$ and $s = \operatorname*{Max}_{1 \leq i \leq m}(s_i)$ is the LFR degree of (6). The polytope $\Delta = \{\Delta(\theta): \theta \in \Theta\}$ has v vertices $\Delta^{(i)}$; $i = \overline{1, v}$, $v \leq \overline{v} = \prod_{i=0}^{r}[v_i]$, which depends on that of the polytope $\Theta$ which parametrizes the system, and the parameter vector $\theta \in \Theta$ parametrizes the whole dynamics of $A(\theta(t))$. An equivalent description to (5)-(6) is

$A(\theta(t)) = \sum_{i=0}^{r} \left( A_{0i} + B_{qi} \Delta(\theta(t))(I - D_{pq}\Delta(\theta(t)))^{-1} C_p U(t - h_i) \right)$ which leads to the multi-input single-output LFR:

$$\dot{x}(t) = \sum_{i=0}^{r} \left( A_{0i} x(t-h_i) + B_{qi} q(t - h_i) \right) \quad ; \quad q(t) = \Delta(\theta) p(t)$$

$$p(t) = C_p x(t) + D_{pq} q(t)$$

$$\Delta(\theta) = \operatorname{Diag}(\theta_1 I_{s_1}, \theta_2 I_{s_2}, ..., \theta_m I_{s_m}) \quad (7)$$

where $q \in \mathbf{R}^d$, $p \in \mathbf{R}^d$, $B_{qi} \in \mathbf{R}^{n \times d}$, $C_p \in \mathbf{R}^{d \times n}$ $(i = \overline{0, r})$. The LFR´s (6)-(7) may be rewritten in the form (4) by defining real matrices:

$$B_q = [B_{q0}, B_{q1}, ..., B_{qr}] \quad ; \quad C_p = [C_{p0}^T, C_{p1}^T, ..., C_{pr}^T]^T$$
$$D_{pq} = \operatorname{Block\ Diag}[D_{pq0}, D_{pq1}, ..., D_{pqr}]$$

$$\Delta(\theta) = \text{Diag}\left[\Delta_0(\theta), \Delta_1(\theta), \ldots, \Delta_r(\theta)\right] \tag{8.a}$$

and

$$p(t) = \left[p_0^T(t), p_1^T(t), \ldots, p_r^T(t)\right]^T$$
$$q(t) = \left[q_0^T(t), q_1^T(t), \ldots, q_r^T(t)\right]^T \tag{8.b}$$

to yield :

$$A(\theta(t)) = \sum_{i=0}^{r}\left(A_{0i} + B_{qi}\Delta_i(\theta(t))\left(I_{d_i} - D_{pqi}\Delta_i(\theta(t))\right)^{-1}C_{pi}U(t-h_i)\right)$$

**Remarks**. (**1**) The LFR´s (6) and (7) are equivalent, the first one describing the delayed dynamics through equivalent output delays while the second one describes it through equivalent input delays.

(**2**) Eqns. 8 prove that the second LFR may be equivalently rewritten in the form (4), associated with the dynamics representations (3). The equivalence arises from the fact that either each $A_i$ may be parametrized by a particular parameter vector $\theta_i \in \Theta_i$ with the matrix function $A(\theta(t)) = \sum_{i=0}^{r} A_i(\theta(t))U(t-h_i)$ being parametrized by

$\theta = (\theta_0^T, \theta_1^T, \ldots, \theta_r^T)^T \in \Theta = \Theta_0 \times \Theta_1 \times \ldots \times \Theta_r$

Since the second LFR may be rewritten as the first one, the formalism presented in this manuscript will be developed for the first LFR with no loss in generality.

(**3**) Note that the overall number of distinct vertices of polytopes v to be checked is $v \leq \bar{v} = \prod_{i=0}^{r}[v_i]$. That inequality may be strict depending on the problem at hand since, depending on the parametrization, some of the matrices defining (1) may eventually generate common vertices. A simple case implying $v \neq \bar{v}$ is when two of those matrices are identical so that the associated set of vertices $\Delta_{(.)}^{(.)}$ become identical. To simplify the notation, we consider in the following $v = \bar{v}$ with no loss in generality noting that the stability for common vertices only require to be tested once.

Note that if the unforced system is globally asymptotically stable independent of delays then $A_{00}$ and $\sum_{i=0}^{r} A_{0i}$ are both stability matrices since $\Theta$ contains the origin $\theta=0$ (describing the uncertainty-free system) and such a system has to be stable for $h_i \to \infty$ and $h_i = 0$ $(i = \overline{0,r})$.

## 3. MAIN RESULTS

The following stability result is concerned with the first LFR eqns. 4. Its proof is omitted.

**Theorem 1**. The (unforced) system (2) is globally asymptotically stable independent of the delays if there exist real matrices $P = P^T > 0$, $S_i = S_i^T > 0$ $(i = \overline{1,r})$ such that, for every $\Delta_i(\theta) \in \Delta_i$, $\theta \in \Theta$ $(i = \overline{0,r})$, there exist complex matrices $G_{\Delta i}(\Delta(\theta)) \in \mathbb{C}^{n \times d_i}$ and $H_{\Delta i}(\Delta(\theta)) \in \mathbb{C}^{d_i \times d_i}$ such that the square $n(r+1)+d$ $\left(d = \sum_{i=0}^{r} d_i\right)$ real symmetric matrix $Q(\Delta(\theta))$ defined as $Q(\Delta(\theta))$ = Block Matrix $\left[Q_{ij}(\Delta(\theta)); i,j=\overline{1,3}\right] < 0$, where

$$Q_{11} = A_{00}^T P + P A_{00} + \sum_{i=1}^{r} S_i + G_{\Delta 0} C_{p0} + C_{p0}^T G_{\Delta 0}^*$$

$$Q_{12} = Q_{21}^T = [PA_{01}, PA_{02}, ..., PA_{0r}]$$

$$Q_{13} = Q_{31}^T$$

$$= \left[P(B_{q0}\Delta_0) + G_{\Delta 0}(D_{pq0}\Delta_0) - G_{\Delta 0} + C_{p0}^T H_{\Delta 0}^* \vdots P(B_{q1}\Delta_1), ..., P(B_{qr}\Delta_r)\right]$$

$$Q_{22} = \text{Block Diag}\left[G_{\Delta 1} C_{p1} + C_{p1}^T G_{\Delta 1}^* - S_1 \vdots ... \vdots G_{\Delta r} C_{pr} + C_{pr}^T G_{\Delta r}^* - S_r\right]$$

$$Q_{23} = Q_{32}^T = \text{Block Diag}$$

$$\left[G_{\Delta 1}(D_{pq1}\Delta_1) - G_{\Delta 1} + C_{p1}^T H_{\Delta 1}^* \vdots ... \vdots G_{\Delta r}(D_{pqr}\Delta_r) - G_{\Delta r} + C_{pr}^T H_{\Delta r}^*\right]$$

$$Q_{33} = \text{Block Diag}\left[H_{\Delta 0}(D_{pq0}\Delta_0) + (D_{pq0}\Delta_0)^T H_{\Delta 0}^* - H_{\Delta 0} - H_{\Delta 0}^* \right.$$

$$\left. \vdots ... \vdots H_{\Delta r}(D_{pqr}\Delta_r) + (D_{pqr}\Delta_r)^T H_{\Delta r}^* - H_{\Delta r} - H_{\Delta r}^*\right]$$

(9)

If $G_{\Delta(.)}$ and $H_{\Delta(.)}$ are restricted to have special forms such that $Q(\Delta(\theta))$ is convex for $\Theta$ being convex, then it suffices that $Q(\Delta(\theta)) < 0$ for its evaluation at each particular j-th vertex $\Delta_i^{(k)}$ $(k=\overline{0,v_i}; i=\overline{0,r})$ generated from the vertices of $\Theta$. Also, if the uncertainty parameter vector is $\theta = 0$ then a well-known simplified version of Theorem 1 for $\Delta_i(0)=0$; $i=\overline{0,r}$, from the last equation in (4), guarantees global asymptotic stability independent of the delays of the (uncertainty-free) nominal system (2). From simple matrix rank continuity arguments, that property is still guaranteed in a robustness stability context within a certain open neighborhood of $\theta = 0$. The subsequent Corollaries to Theorem 1 follow.

**Corollary 1**. The (unforced) system (2) is globally asymptotically stable independent of the delays if there exist real matrices $P = P^T > 0$, $S_i = S_i^T > 0$ $(i=\overline{1,r})$ such that, for every $\Delta_i(\theta) \in \Delta_i$, $\theta \in \Theta$ $(i=\overline{0,r})$, there exist complex matrices $G_{\Delta i}(\Delta(\theta)) \in C^{n \times d_i}$ and $H_{\Delta i}(\Delta(\theta)) \in C^{d_i \times d_i}$ such that the square $n(r+1)$ real symmetric matrix $Q'(\theta)$ defined as $Q'(\Delta(\theta)) = \text{Block Matrix}\left[Q'_{ij}(\Delta(\theta)); i,j=1,2\right] < 0$, where

$$Q'_{11} = A_{00}^T P + P A_{00} + \sum_{i=1}^{r} S_i + \tilde{Q}'_{11}$$

$$Q'_{12} = Q'_{21}{}^T = [PA_{01}, PA_{02}, ..., PA_{0r}] + \tilde{Q}'_{12}$$

$$Q'_{22} = \text{Block Diag}[-S_1, ..., -S_r] + \tilde{Q}'_{22}$$

(10.a)

where

$$\tilde{Q}'_{11} = G_{\Delta 0} C_{p0} + C_{p0}^T G_{\Delta 0}^* + \left[P(B_{q0}\Delta_0)\right.$$

$$\left. + G_{\Delta 0}(D_{pq0}\Delta_0) - G_{\Delta 0} + C_{p0}^T H_{\Delta 0}^*\right]$$

$$\left(I - D_{pq0} \Delta_0\right)^{-1} C_{p0} + C_{p0}^T \left(I - D_{pq0} \Delta_0\right)^{-T}$$
$$\left[\left(B_{q0} \Delta_0\right)^T P + \left(D_{pq0} \Delta_0\right)^T G_{\Delta 0}^* - G_{\Delta 0}^* + H_{\Delta 0} C_{p0}\right]$$

$$\tilde{Q}'_{12} = \tilde{Q}'_{21}{}^T =$$

$$\left[P\left(B_{q1} \Delta_1\right)\left(I - D_{pq1} \Delta_1\right)^{-1} C_{p1} \vdots \ldots \vdots P\left(B_{qr} \Delta_r\right)\left(I - D_{pqr} \Delta_r\right)^{-1} C_{pr}\right]$$

$$\tilde{Q}'_{22}{}^{(i)} = G_{\Delta i} C_{pi} + C_{pi}^T G_{\Delta i}^* + \left[G_{\Delta i}\left(D_{pqi} \Delta_i\right) - G_{\Delta i} + C_{pi}^T H_{\Delta i}^*\right]$$
$$\left(I - D_{pqi} \Delta_i\right)^{-1} C_{pi}$$

$$+ C_{pi}^T \left(I - D_{pqi} \Delta_i\right)^{-T} \left[\left(D_{pqi} \Delta_i\right)^T G_{\Delta i}^* - G_{\Delta i}^* + H_{\Delta i}^* C_{pi}\right] +$$
$$C_{pi}^T \left(I - D_{pqi} \Delta_i\right)^{-T} \left[H_{\Delta i}\left(D_{pqi} \Delta_i\right) + \left(D_{pqi} \Delta_i\right)^T H_{\Delta i}^* - H_{\Delta i} - H_{\Delta i}^*\right]$$
$$\left(I - D_{pqi} \Delta_i\right)^{-1} C_{pi} \quad \left(i = \overline{1, r}\right) \tag{10.b}$$

**Corollary 2**. Assume that the system (4) is time-invariant, nominally parametrized at $\theta^* = 0 \in \Theta$ with no parametrical uncertainties; i.e. $\theta(t) = \theta^*(t) = 0$ for all time. Then, it is globally asymptotically stable independent of the delays if the real symmetric $(r+1)n$ – square matrix

$$Q_0^{(2)} = \text{Block Matrix}\left[Q_{0ij}^{(2)}; i, j = 1, 2\right] < 0 \text{ where}$$

$$Q_{011}^{(2)} = A_{00}^T P + P A_{00} + \sum_{i=1}^{r} S_i$$

$$Q_{012}^{(2)} = Q_{021}^{(2)}{}^T = Q_{12} = \left[P A_{01}, P A_{02}, \ldots, P A_{0r}\right]$$

$$Q_{022}^{(2)} = \text{Block Diag}\left[-S_1, \ldots, -S_r\right]$$

(11)

for some real $P = P^T > 0$; $S_i = S_i^T > 0$ $\left(i = \overline{1, r}\right)$. A necessary condition for $Q_0^{(2)} < 0$ is that $A_{00}$ be a stability matrix (i.e. of all its eigenvalues of negative real parts). Also, for some real $R > 0$, there exists an open neigborhood $N(0, R)$ of $0 \in \Theta$ (the polytope of parametrical uncertainties) of radius R such that the system (4) ( and then the unforced system (2) ) is globally asymptotically stable independent of the delays if $\theta \in \{N(0, R) \cap \Theta\}$ for all time $t \geq 0$.

**Corollary 3**. The (unforced) system (2) is globally asymptotically stable independent of the delays if there exist if there exist $v = \prod_{i=0}^{r}\left[v_i\right]$ matrices

$$Q^{(3)}\left(k_0, k_1, \ldots, k_r\right) = \text{Block Matrix}\left[Q_{ij}^{(3)}\left(k_0, k_1, \ldots, k_r\right); i, j = \overline{1, 3}\right] < 0$$

for all $k_i = \overline{1, v_i}$; $i = \overline{0, r}$ for some real n-matrices $P = P^T > 0$, $M_i = M_i^T > 0$; $i = \overline{0, r}$ where

$$Q_{11}^{(3)} = A_{00}^T P + P A_{00} + \sum_{i=1}^{r} S_i + C_{p0}^T M_0 C_{p0}$$

$$Q_{12}^{(3)} = Q_{21}^{(3)}{}^T = Q_{12} = \left[P A_{01}, P A_{02}, \ldots, P A_{0r}\right]$$

$$Q_{13}^{(3)}(k_0,...,k_r) = Q_{31}^{(3)T}(k_0,...,k_r)$$

$$= \left[ P\left(B_{q0}\Delta_0^{(k_0)}\right) + C_{p0}^T M_0\left(D_{pq0}\Delta_0^{(k_0)}\right) \vdots \right.$$

$$\left. P\left(B_{q1}\Delta_1^{(k_1)}\right),...,P\left(B_{qr}\Delta_r^{(k_r)}\right) \right]$$

$$Q_{22}^{(3)} = \text{Block Diag}\left[ C_{p1}^T M_1 C_{p1} - S_1,...,C_{pr}^T M_r C_{pr} - S_r \right]$$

$$Q_{23}^{(3)}(k_1,...,k_r) = Q_{32}^{(3)}(k_1,...,k_r)^T$$

$$= \text{Block Diag}\left[ C_{p1}^T M_1\left(D_{pq1}\Delta_1^{(k_1)}\right),...,C_{pr}^T M_r\left(D_{pqr}\Delta_r^{(k_r)}\right) \right]$$

$$Q_{33}^{(3)}(k_0,...,k_r) = \text{Block Diag}$$

$$\left[ -M_0 + \left(D_{pq0}\Delta_0^{(k_0)}\right)^T M_0\left(D_{pq0}\Delta_0^{(k_0)}\right) \vdots ... \vdots \right.$$

$$\left. -M_r + \left(D_{pqr}\Delta_r^{(k_r)}\right)^T M_r\left(D_{pqr}\Delta_r^{(k_r)}\right) \right] \quad (12)$$

**Remark 4**. Note that Corollary 3 may be tested for any set of real symmetric positive definite matrices $M_i$ $(i = \overline{0,r})$ and holds, in particular, if the stability test becomes positive for the v-matrices $Q^{(3)}(k_0,...,k_r) < 0$ for some identical symmetric matrices $M_i = M > 0$ $(i = \overline{0,r})$. Corollary 3 adopts the following parallel form, under weaker conditions, if the LFR degrees are unity with respect to all the delays.

**Corollary 4**. Assume that $s_i = 1$; $i = \overline{0,r}$. Then, the unforced system (4) is globally asymptotically stable independent of the delays if there exist $P = P^T > 0$ and v real positive definite symmetric matrices $M(k_i)$; $k_i = \overline{1,v_i}$; $i = \overline{0,r}$ such that the v real symmetric matrices :

$$Q^{(4)}(k_0,k_1,..., k_r) = \text{Block Matrix}\left[Q_{ij}^{(4)}(k_0,k_1,..., k_r); i,j = \overline{1,3}\right] < 0$$

where the block matrices have the same structures as in Corollary 3 with the replacements $M_i \to M(k_i)$; $k_i = \overline{1,v_i}$; $i = \overline{0,r}$.

It is obvious that Corollary 4 is stronger than Corollary 3 since the v matrix negative definiteness tests might be tested with $M(k_i)$ distinct matrices (one per vertex). Parallel results were first obtained for the delay-free case in [3].

## 4. ASYMPTOTIC STABILITY DEPENDENT ON THE DELAYS

Parallel results to those in the above Section II may be obtained depending on the first intervals for delays ensuring global asymptotic stability; i.e $h_i \in \left[0, h_i^0\right]$ $(i = \overline{1,r})$ and $h_0 = h_0^0 = 0$. The results are obtained from Lyapunov's function candidates that include integral terms for the effects of coupled combined delays. The ´ad-hoc´ version of Theorem 1 for this situation is:

**Theorem 2**. The (unforced) system (2) is globally asymptotically stable for all the delays $h_i \in [0, h_i^0]$ ($i = \overline{1,r}$) if there exist real matrices $P = P^T > 0$, $S_{ij} = S_{ij}^T > 0$ ($i = \overline{1,r}$; $j = \overline{0,r}$) such that, for every $\Delta_i(\theta) \in \Delta_i$, $\theta \in \Theta$ ($i = \overline{0,r}$), there exist complex matrices $G_{\Delta ij}(\Delta(\theta)) \in C^{n \times d_i}$ and $H_{\Delta ij}(\Delta(\theta)) \in C^{d_i \times d_i}$ such that the square $n(r+1)r + d$ $\left(d = \sum_{i=0}^{r} d_i\right)$ real symmetric matrix $\hat{Q}(\Delta(\theta))$ defined as $\hat{Q}(\Delta(\theta)) =$ Block Matrix

$\left[\hat{Q}_{ij}(\Delta(\theta)) ; i,j = \overline{1,3}\right] < 0$, where

$$\hat{Q}_{ij} = \hat{Q}_{0ij} + \tilde{\hat{Q}}_{0ij} \quad ; \quad i, j = \overline{0, r}$$

$$\hat{Q}_{011} = \left(\sum_{i=0}^{r} A_{0i}^T\right) P + P\left(\sum_{i=0}^{r} A_{0i}\right) + \sum_{i=1}^{r}\sum_{j=0}^{r} h_i^0 S_{ij}$$

$$\hat{Q}_{012} = \hat{Q}_{021}^T = \left[h_1^0 P A_{01} \hat{A}, ..., h_r^0 P A_{0r} \hat{A}\right]$$

$$\hat{Q}_{013} = \hat{Q}_{031}^T = \left[P(B_{q0} \Delta_0), P(B_{q1} \Delta_1), ..., P(B_{qr} \Delta_r)\right]$$

$$\hat{Q}_{022} = \text{Block Diag}\left[-h_1^0 R_1, ..., -h_r^0 R_r\right] \quad (14)$$

$$\hat{Q}_{023} = \hat{Q}_{032}^T = 0 \quad ; \quad \tilde{\hat{Q}}_{033} = 0$$

$$\hat{A} = [A_0, A_1, ..., A_r] \in \mathbf{R}^{n \times n(r+1)}$$

$$R_i = [S_{i0}, S_{i1}, ..., S_{ir}] \in \mathbf{R}^{n(r+1) \times n(r+1)}, \text{ for } i = \overline{1, r} \quad (15)$$

$$\hat{Q}_{11} = \hat{Q}_{11}^0 \quad ; \quad \hat{Q}_{12} = \hat{Q}_{012}^0$$

$$\hat{Q}_{13} = \hat{Q}_{31}^T = \Big[P(B_{q0} \Delta_0) + G_{\Delta 0}(D_{pq0} \Delta_0) - G_{\Delta 0} + C_{p0}^T H_{\Delta 00}^* \vdots$$

$$\underbrace{P(B_{q1} \Delta_1)}_{nq}, \underbrace{0, ..., 0}_{rnq} \vdots ... \vdots \underbrace{P(B_{qr} \Delta_r)}_{nq}, \underbrace{0, ..., 0}_{rnq}\Big]$$

$$\hat{Q}_{23} = \hat{Q}_{32}^T = \text{Block Diag }\Big[G_{\Delta 10}(D_{pq1} \Delta_1) - G_{\Delta 10} + C_{p1}^T H_{\Delta 10}^* \vdots ... \vdots$$

$$G_{\Delta 1r}(D_{pq1} \Delta_1) - G_{\Delta 1r} + C_{p1}^T H_{\Delta 1r}^* \vdots ... \vdots$$

$$G_{\Delta r0}(D_{pqr} \Delta_r) - G_{\Delta r0} + C_{pr}^T H_{\Delta r0}^* \vdots ... \vdots$$

$$G_{\Delta rr}(D_{pqr} \Delta_r) - G_{\Delta rr} + C_{pr}^T H_{\Delta rr}^*\Big]$$

$$\hat{Q}_{33} = \text{Block Diag }\Big[H_{\Delta 10}(D_{pq1} \Delta_1) + (D_{pq1} \Delta_1)^T H_{\Delta 10}^*$$

$$-H_{\Delta 10} - H_{\Delta 10}^* \vdots ... \vdots H_{\Delta 1r}(D_{pq1} \Delta_1) + (D_{pq1} \Delta_1)^T H_{\Delta 1r}^* - H_{\Delta 1r} - H_{\Delta 1r}^* \vdots ... \vdots$$

$$H_{\Delta r0}(D_{pqr} \Delta r) + (D_{pqr} \Delta_r)^T H_{\Delta r0}^* - H_{\Delta r0} - H_{\Delta r0}^* \vdots ... \vdots$$

$$H_{\Delta rr}(D_{pqr} \Delta_r) + (D_{pqr} \Delta_r)^T H_{\Delta rr}^* - H_{\Delta rr} - H_{\Delta rr}^*\Big]$$

(16)

Note that, since $\hat{Q}_{11} < 0$ then the delay-free system, $\dot{z}(t) = \left(\sum_{i=0}^{r} A_{0i}\right) z(t)$ is globally exponentially stable if Theorem 2 holds. Since (4) is a LFR of (2), the global asymptotic stability of (3), dependent on the delays, may be analyzed for all (constant or time-varying) $\theta \in \Theta$ from that of the LFR (4) by considering the vertices of the polytope $\Delta$. Therefore, the following Corollaries 5-6 to Theorem 2 may be stated and proved in a similar way as Corollaries 3-4 to Theorem 1.

**Corollary 5**. The (unforced) system (2) is globally asymptotically stable independent of the delays if there exist if there exist (at most) $v = \prod_{i=0}^{r} [v_i]$ matrices

$\hat{Q}^{(5)}(k_0, k_1, \ldots, k_r) = \text{Block Matrix} \left[ \hat{Q}_{ij}^{(5)}(k_0, k_1, \ldots, k_r); i,j = \overline{1,3} \right] < 0$

for all $k_i = \overline{1, v_i}$; $i = \overline{0, r}$ for all the delays $h_i \in [0, h_i^0]$ $(i = \overline{0,r})$ with $h_0 = h_0^0 = 0$ and some real n-matrices $P = P^T > 0$, $M_i = M_i^T > 0$; $i = \overline{0,r}$, where:

$$\hat{Q}_{11}^{(5)} = \left( \sum_{i=0}^{r} A_{0i}^T \right) P + P \left( \sum_{i=0}^{r} A_{0i} \right) + \sum_{i=1}^{r} \sum_{j=0}^{r} h_i^0 S_{ij} + C_{p0}^T M_0 C_{p0}$$

$$\hat{Q}_{12}^{(5)} = \hat{Q}_{21}^{(5)T} = Q_{12}^{(3)}; \quad \hat{Q}_{13}^{(5)}(k_0, \ldots, k_r) = \hat{Q}_{31}^{(5)T} = Q_{13}^{(3)}(k_0, \ldots, k_r)$$

(defined in eqns. 12)

$$\hat{Q}_{13}^{(5)}(k_0, \ldots, k_r) = \hat{Q}_{31}^{(5)T} = \left[ P \left( B_{q0} \Delta_0^{(k_0)} \right) + C_{p0}^T M_0 \left( D_{pq0} \Delta_0^{(k_0)} \right), \right.$$

$$\left. P \left( B_{q1} \Delta_1^{(k_1)} \right), \ldots, P \left( B_{qr} \Delta_r^{(k_r)} \right) \right]$$

$$\hat{Q}_{22}^{(5)} = \text{Block Diag} \left[ C_{p1}^T M_1 C_{p1} - h_1^0 S_{10}, -S_{11}, \ldots, -S_{1r} \vdots \ldots \vdots \right.$$

$$\left. C_{pr}^T M_r C_{pr} - h_r^0 S_{r0}, -S_{r1}, \ldots, -S_{rr} \right]$$

$$\hat{Q}_{23}^{(5)}(k_1, \ldots, k_r) = \hat{Q}_{32}^{(5)}$$

$$= \text{Block Diag} \left[ C_{p1}^T M_1 \left( D_{pq1} \Delta_1^{(k_1)} \right) \vdots \underbrace{0}_{rn \times rn} \vdots \ldots \vdots C_{pr}^T M_r \left( D_{pqr} \Delta_r^{(k_r)} \right) \vdots \underbrace{0}_{rn \times rn} \right]$$

$$\hat{Q}_{33}^{(5)}(k_1, \ldots, k_r) = Q_{33}^{(3)}(k_1, \ldots, k_r) \tag{17}$$

ACKNOWLEDGMENT
The author is grateful to MCYT by its partial support of this work trough Grant DPI 2006-00714